\documentclass{amsart}
\usepackage{setspace}

\usepackage{amsmath}
\usepackage{amsthm}
\usepackage{amssymb}

\usepackage[latin1]{inputenc} 

\newtheorem{theorem}{Theorem}[section]

\newtheorem{proposition}[theorem]{Proposition}
\newtheorem{corollary}[theorem]{Corollary}
\newtheorem{remark}[theorem]{Remark}

\newtheorem{definition}[theorem]{Definition}
\newtheorem{example}[theorem]{Example}

\numberwithin{equation}{section}

\begin{document}
\title[A type of the entropy of an ideal]
{A type of the entropy of an ideal}

\author{Nicu\c{s}or Minculete}
\address{Faculty of Mathematics and Computer Science, Transilvania University\\
 Iuliu Maniu street 50, Bra\c{s}ov 500091, Romania}
\email{minculete.nicusor@unitbv.ro}

\author{Diana Savin}
\address{Faculty of Mathematics and Computer Science, Transilvania University\\
 Iuliu Maniu street 50, Bra\c{s}ov 500091, Romania}
\email{diana.savin@unitbv.ro; dianet72@yahoo.com}

\date{}
\begin{abstract}
In this article we find some properties of certain types of entropies of a natural number. Also, regarding the entropy $H$ of a natural number, introduced by Minculete and Pozna,
we generalize this notion for ideals and we find some of its properties. In the last section we find some inequalities, involving the  entropy  $H$ of an  exponential divisor of a positive integer, respectively 
the  entropy  $H$ of an  exponential divisor of an ideal.
\end{abstract}
\maketitle
\section{Introduction and Preliminaries} 
\noindent 
\smallskip\\

In information theory, the entropy is defined as a measure of uncertainty. The most used of entropies is the Shannon entropy ($H_S$), which is given for a probability distribution ${\bf p}=\{p_1,...,p_r\}$ thus
$$H_S({\bf p})= - \sum_{i=1}^{r}p_i\cdot \log p_i.$$
We have found several ways to define the entropy of a natural number. Jeong {\it et al.}, in \cite{Jeong}, defined the additive entropy of a natural number in terms of the additive partition function. In \cite{link}, we found the following definition for the entropy of a natural number:
$$\overline{H}(n):= - \sum_{d|n}^{}\frac{d}{\sigma(n)}\log \frac{d}{\sigma(n)}=\log\sigma(n)-\frac{1}{\sigma(n)}\sum_{d|n}^{}d\log d,$$
where $\sigma(n)$ is the sum of natural divisors of $n$. This entropy has the following interesting property: $$\overline{H}(mn)=\overline{H}(m)+\overline{H}(n),$$
when $m,n\in\mathbb{N}^*$ and $\gcd(m,n)=1$. If $p$ is a prime number and $\alpha\in\mathbb{N}^*$, then we have $$\overline{H}(p^\alpha)=-\frac{(\alpha+1)\log p}{p^{\alpha+1}-1}+\log\frac{1-p^{-(\alpha+1)}}{p-1}+\frac{p\log p}{p-1}.$$
Taking the limit as $\alpha\to\infty$, we obtain
\begin{equation}\label{lim_1}
\lim_{\alpha\to\infty}\overline{H}(p^\alpha)=\frac{p\log p}{p-1}-\log(p-1).
\end{equation}
We remark that, if $p$ is a prime number, $q>1$ such that $\frac{1}{p}+\frac{1}{q}=1$, then
$$H_S\left(\frac{1}{p},\frac{1}{q}\right)=\frac{p-1}{p}\left(\frac{p\log p}{p-1}-\log(p-1)\right)=\left(1-\frac{1}{p}\right)\lim_{\alpha\to\infty}\overline{H}(p^\alpha).$$
In the paper \cite{Minculete1}, Minculete and  Pozna introduced the notion of entropy of a natural number by another way, namely: if $n\in\mathbb{N},$  $n\geq2,$ applying the Fundamental Theorem of Arithmetic, $n$ is written uniquely
 $n=p^{\alpha_{1}}_{1}p^{\alpha_{2}}_{2}...p^{\alpha_{r}}_{r}$, where $r\in\mathbb{N}^{*}$, $p_{1}, p_{2},...,p_{r}$ are distinct prime positive integers and $\alpha_{1}, \alpha_{2},...,\alpha_{r}\in\mathbb{N}^{*}.$ Let $\Omega\left(n\right)=\alpha_{1}+\alpha_{2}+...+\alpha_{r}$ and $p\left(\alpha_{i}\right)=\frac{\alpha_{i}}{\Omega\left(n\right)},$
$\left(\forall\right)$ $i=\overline{1,r}$. The entropy of $n$  is defined by
\begin{equation}
H\left(n\right)= - \sum_{i=1}^{r}p\left(\alpha_{i}\right)\cdot \log \: p\left(\alpha_{i}\right),    \tag{1.1}
\end{equation}
where $\log$ is the natural logarithm. Here, by convention $H(1)=0$.\\
Minculete and  Pozna  (in  \cite{Minculete1}) gave an equivalent form for the entropy of $n,$ namely:
\begin{equation}
H\left(n\right)=\log \: \Omega\left(n\right) - \frac{1}{\Omega\left(n\right)}\cdot \sum_{i=1}^{r}\alpha_{i}\cdot \log \: \alpha_{i}.  \tag{1.2}
\end{equation}
For example, if $n=6=2\cdot 3,$ we have:
$$ H\left(6\right)=\log \: 2 -  \frac{1}{2}\cdot 2\cdot \log \: 1=\log \: 2=0.6931.... $$
Another example: if $n=24=2^{3}\cdot 3,$ we have:
$$ H\left(24\right)=\log \: 4 -  \frac{1}{4}\cdot 3\cdot \log \: 3= \frac{1}{4}\cdot \log\left(\frac{4^4}{3^3}\right)=2.2493.... $$
Minculete and  Pozna proved (in  \cite{Minculete1}) the followings:
\begin{proposition}
\label{onedotone}
\begin{equation}
0\leq H\left(n\right)\leq \log \: \omega\left(n\right), \; \left(\forall\right) \; n\in \mathbb{N},\; n\geq 2,  \tag{1.3}
\end{equation}
where $ \omega\left(n\right)$  is the number of distinct prime factors of $n.$
\end{proposition}
\begin{remark}
\label{onedottwo}
i) If $n =p^{\alpha}$ , then $H\left(n\right)=0;$\\
ii) If $n =p_{1}\cdot p_{2}\cdot...\cdot p_{r},$ then $H\left(n\right)=log\:\omega\left(n\right);$\\
iii) If $n =\left(p_{1}\cdot p_{2}\cdot...\cdot p_{r}\right)^{k},$ then $H\left(n\right)=log\:\omega\left(n\right).$
\end{remark}
It is easy to see that $H(n^{\alpha})=H(n)$, with $\alpha\geq 1$.\\
\section{A comparison between the entropies $H$ and $\overline{H}$} 
\noindent 
\smallskip\\
In this section we propose to compare the entropies $H$ and $\overline{H}$ looking to similarities and differences between them.
\begin{proposition}
\label{ms}
\begin{equation}\label{sec2_1}
\lim_{p\to\infty}\lim_{\alpha\to\infty}\overline{H}(p^\alpha)=0. 
\end{equation}
\end{proposition}
\begin{proof}
From relation \eqref{lim_1} we have $\lim_{\alpha\to\infty}\overline{H}(p^\alpha)=\frac{p\log p}{p-1}-\log(p-1).$ Next, we use the following limit of functions:
$$\lim_{x\to\infty}\left(\frac{x\log x}{x-1}-\log(x-1)\right)=\lim_{x\to\infty}\frac{x\log x-(x-1)\log (x-1)}{x-1}$$
$$=\lim_{x\to\infty}\left(\log x-\log (x-1)\right)=\lim_{x\to\infty}\log\frac{x}{x-1}=0.$$
Therefore, we get $\lim_{p\to\infty}\lim_{\alpha\to\infty}\overline{H}(p^\alpha)=\lim_{p\to\infty}\left(\frac{p\log p}{p-1}-\log(p-1)\right)=0.$
\end{proof}
\begin{remark}
It is easy to see that $\lim_{p\to\infty}\lim_{\alpha\to\infty}\overline{H}(p^\alpha)=0=H(p^\alpha).$ 
\end{remark}
\begin{proposition}
\label{ms.1}
If $\gcd(n,p)=1$, with $p$ is a prime number and $n,\alpha\in\mathbb{N}^*$, then we have 
\begin{equation}\label{sec2_2}
\lim_{\alpha\to\infty}H(np^\alpha)=0. 
\end{equation}
\end{proposition}
\begin{proof}
From the definition of $H$ we have 
\begin{align*}H(np^\alpha)&=\log(\Omega(n)+\alpha)-\frac{1}{\Omega(n)+\alpha}\left(\sum_{i=1}^{r}\alpha_{i}\cdot \log \: \alpha_{i}+\alpha\log\alpha\right)\\
&=\log(\Omega(n)+\alpha)-\frac{\alpha\log\alpha}{\Omega(n)+\alpha}-\frac{1}{\Omega(n)+\alpha}\left(\Omega(n)\log\Omega(n)-\Omega(n)H(n)\right)\\
&=\frac{\Omega(n)H(n)}{\Omega(n)+\alpha}+\log(\Omega(n)+\alpha)-\frac{\Omega(n)\log \Omega(n)+\alpha\log\alpha}{\Omega(n)+\alpha}.
\end{align*}
It follows that 
\begin{equation}\label{ms_5}
H(np^\alpha)=\frac{\Omega(n)H(n)}{\Omega(n)+\alpha}+\log(\Omega(n)+\alpha)-\frac{\Omega(n)\log \Omega(n)+\alpha\log\alpha}{\Omega(n)+\alpha}.
\end{equation}
Taking to limit when $\alpha\to\infty$, we deduce the relation of the statement.
\end{proof}
\begin{remark}
Related to the entropy $\overline{H}$ we have $$\lim_{\alpha\to\infty}\overline{H}(np^\alpha)=\overline{H}(n)+\frac{p\log p}{p-1}-\log(p-1),$$ 
when $\gcd(n,p)=1$, with $p$ is a prime number and $n,\alpha\in\mathbb{N}^*$.
\end{remark}
We also see that if $\gcd(m,n)=1$, then $$H(mn)\neq H(m)+H(n).$$ As a result, we ask ourselves the question of what is the relationship between $H(mn)$ and $H(m)+H(n)$, where $m,n\in\mathbb{N}^*$, $m,n\geq 2$.

If $m=22$ and $n=105$, then $H(m)=\log 2$, $H(n)=\log 3$ and $H(mn)=\log 5$, so, we have
$$H(mn)< H(m)+H(n).$$

If $m=20$ and $n=63$, then $H(m)=H(n)=\log 3-\frac{2}{3}\log 2$ and $H(mn)=\log 6-\frac{2}{3}\log 2$, which means that $$H(mn)-H(m)-H(n)=\frac{1}{3}\left(5\log 2-3\log 3\right)=\frac{1}{3}\log\frac{32}{27}>0$$
so, we have
$$H(mn)>H(m)+H(n).$$

Next, we study a general result of this type for the entropy $H$.
\begin{proposition}
\label{ms.2}
We assume that $m=p^kq$ and $n=p^kt$, where $p,q,t$ are distinct prime numbers and $k\in\mathbb{N}^*$. Then the inequality 
$$H(mn)<H(m)+H(n)$$
holds.
\end{proposition}
\begin{proof}
From the definition of $H$ we have $H(m)=H(n)=\log (k+1)-\frac{k}{k+1}\log k$ and $H(mn)=\log 2(k+1)-\frac{k}{k+1}\log 2k$. Therefore, we obtain
$$H(m)+H(n)-H(mn)=\frac{1}{k+1}\left((k+1)\log (k+1)-k\log k-\log 2\right).$$
We consider the function $f:[1,\infty)\to\mathbb{R}$ defined by $f(x)=(x+1)\log (x+1)-x\log x-\log 2$. Since $f'(x)=\log\frac{x+1}{x}>0$ for every $x\geq 1$, we deduce that
the function $f$ is increasing, so, we have $f(x)\geq f(1)=\log 2>0.$
Consequently, the inequality of the statement is true.
\end{proof} 
\begin{proposition}
\label{ms.3}
We assume that $m=p_1^kp_2$ and $n=q_1^kq_2$, where $p_1,p_2,q_1,q_2$ are distinct prime numbers and $k\in\mathbb{N}^*$. Then we have the following inequality 
$$H(mn)\geq H(m)+H(n).$$
Equality holds for $k=1$.
\end{proposition}
\begin{proof}
For $k=1$, we deduce that $m=p_1p_2$ and $n=q_1q_2$, which implies $H(m)=H(n)=\log 2$ and $H(mn)=\log 4$, so, we have $$H(mn)=H(m)+H(n).$$
For $k\geq 2$, we find $H(m)=H(n)=\log (k+1)-\frac{k}{k+1}\log k$ and $H(mn)=\log 2(k+1)-\frac{k}{k+1}\log k$. Now, we obtain
$$H(mn)-H(m)-H(n)=\frac{1}{k+1}\left((k+1)\log 2 +k\log k-(k+1)\log (k+1)\right)$$
for all $k\geq 2$, because the function $f:[2,\infty)\to\mathbb{R}$ defined by $f(x)=(x+1)\log 2 +x\log x-(x+1)\log (x+1)$ is strictly positive. It is easy to see that $f'(x)>0$ for every $x\geq 2$. Therefore, for $x=k$, we prove the relation of the statement.
\end{proof} 
We study another result for which we have $$H(mn)\geq H(m)+H(n),$$
where $m,n\in\mathbb{N}^*$, $m,n\geq 2$.
\begin{proposition}
\label{ms.4}
Let $m,n$ be two natural numbers such that $\gcd(m,n)=1$ and decomposition in prime factors of $m,n$ given by $m=\prod_{i=1}^rp_i^{a_i}$ and $n=\prod_{j=1}^sq_j^{b_j}$ with $a_i,b_j\geq 3$ for all $i\in\{1,...,r\}$ and $j\in\{1,...,s\}$. Then the inequality 
$$H(mn)> H(m)+H(n)$$
holds.
\end{proposition}
\begin{proof}
Using the definition of $H$, we deduce the equality
\begin{equation}\label{eq}
H(mn)-H(m)-H(n)
\end{equation}
$$=\frac{\Omega(n)}{\Omega(m)+\Omega(n)}\sum_{i=1}^ra_i\log a_i+\frac{\Omega(m)}{\Omega(m)+\Omega(n)}\sum_{j=1}^sb_j\log b_j-\log\frac{\Omega(m)\Omega(n)}{\Omega(m)+\Omega(n)}.
$$
Since $\log a_i,\log b_j>1$ for all $i\in\{1,...,r\}$ and $j\in\{1,...,s\}$, we obtain that $\sum_{i=1}^ra_i\log a_i>\sum_{i=1}^ra_i=\Omega(m)$ and $\sum_{j=1}^sb_j\log b_j>\sum_{j=1}^sb_j=\Omega(n)$.
Using equality \eqref{eq} and above inequalities, we show that
$$H(mn)-H(m)-H(n)>\frac{2\Omega(m)\Omega(n)}{\Omega(m)+\Omega(n)}-\log\frac{\Omega(m)\Omega(n)}{\Omega(m)+\Omega(n)}.
$$
We consider the function $f:[\frac{1}{2},\infty)\to\mathbb{R}$ defined by $f(x)=2x-\log x$. Since $f'(x)=\frac{2x-1}{x}\geq 0$ for every $x\geq \frac{1}{2}$, then $f(x)\geq f(\frac{1}{2})=1-\log\frac{1}{2}=1+\log 2>0.$
It is easy to see that $\Omega(m),\Omega(n)\geq 1$. Therefore, we obtain $2\Omega(m)\Omega(n)-\Omega(m)-\Omega(n)=\left(\Omega(m)-1\right)\left(\Omega(n)-1\right)+\Omega(m)\Omega(n)-1\geq 0,$ which means that $\frac{\Omega(m)\Omega(n)}{\Omega(m)+\Omega(n)}\geq \frac{1}{2}.$ So, we have $$\frac{2\Omega(m)\Omega(n)}{\Omega(m)+\Omega(n)}-\log\frac{\Omega(m)\Omega(n)}{\Omega(m)+\Omega(n)}>0.
$$
Consequently, the inequality of the statement is true.
\end{proof} 

Next, our goal was to show that the entropy $H$ is more suitable to extend it to ideals.

\section{The entropy of an ideal}

\noindent 
\smallskip\\
In this section we introduce the notion of entropy of an ideal of a ring of algebraic integers and we find interesting properties of it.

In the $4^{\it rd}$ section we obtain some properties of the entropy of an exponential divisor.

Subbarao  introduced the notion of \textit{exponential divisor} of a positive integer and he found some properties of these divisors (see \cite{Straus}, \cite{Subbarao}). So, if $n$ is a positive integer,  $n>1,$ it can be written uniquely
as $n=p_1^{\alpha_1}p_2^{\alpha_2}\cdot\cdot\cdot p_r^{\alpha_r}$, where $r\in\mathbb{N}^{*},$ $p_{1}, p_{2},...,p_{r}$ are distinct prime positive integers and $\alpha_{1}, \alpha_{2},...,\alpha_{r}\in\mathbb{N}^{*}.$
The positive integer $d=p_1^{\beta_1}p_2^{\beta_2}\cdot\cdot\cdot p_r^{\beta_r}$ (with $\beta_{1}, \beta_{2},...,\beta_{r}\in\mathbb{N}^{*}$) is called \textit{exponential divisor} or \textit{e-divisor} of $n=p_1^{\alpha_1}p_2^{\alpha_2}\cdot\cdot\cdot p_r^{\alpha_r}$, 
if $\beta_i|\alpha_i$, for every $i\in \{1,...,r\}.$
The number of exponential divisors of $n$ is denoted by $\tau^{\left( e\right)}\left(n\right)$ and if $n>1$ denoted as above,  we have $\tau^{\left( e\right)}\left(n\right)=\tau\left( \alpha_1\right)\tau\left( \alpha_2\right)\cdot\cdot\cdot\tau\left( \alpha_r\right),$
where $\tau\left( \alpha_i\right)$ is the number of natural divisors of $ \alpha_i,$  $\left(\forall\right)$ $i=\overline{1,r}.$ By convention $\tau^{\left( e\right)}\left(1\right)=1.$\\
\indent Many properties of the number of the exponential divisors of a positive integer $n$ can be found in the articles \cite{Minculete}, \cite{Min_Sav}, \cite{Pet}, \cite{Sandor_2}, \cite{Sandor}, \cite{Toth}, \cite{Wu}.

Let $K$ be an algebraic number field of degree $[K :
\mathbb{Q}] = n,$ where  $n\in\mathbb{N}$, $n\geq2$ and let $\mathcal{O}_{K}$ be its ring of integers. Let $p$ be a prime positive integer. Since $\mathcal{O}_{K}$ is a Dedekind ring, applying the fundamental theorem of Dedekind rings, 
the ideal $p\mathcal{O}_{K} $ is written uniquely (except for the order of the factors) like this:
$$p\mathcal{O}_{K}=P^{e_{1}}_{1}\cdot P^{e_{2}}_{2}\cdot...\cdot P^{e_{g}}_{g},$$
where  $g\in\mathbb{N}^*,$   $e_{1}, e_{2},..., e_{g}\in\mathbb{N}^*$  and $P_{1},$ $P_{2},$..., $P_{g}\in Spec\left(\mathcal{O}_{K}\right).$
 The number $e_{i}$ ($i = \overline{1,g}$) is called the ramification index of $p$ at the ideal $P_{i}$. 
The following result is known (see \cite{ireland}):

\begin{proposition}
\label{onedotthree}
In the above notations, we have:\\
i)
$$\sum_{i=1}^{g}e_{i}f_{i}=[K :\mathbb{Q}] = n,$$
where  $f_{i}$ is the residual degree of $p,$ meaning $f_{i}=\left[\mathcal{O}_{K}/P_{i} : \mathbb{Z}/p\mathbb{Z}\right],$ $i = \overline{1,g}.$\\ 
ii) If moreover $\mathbb{Q} \subset K$ is a Galois extension, then $e_{1}=e_{2}=...=e_{g}$ (denoted by $e$), $f_{1}=f_{2}=...=f_{g}$ (denoted by $f$).
Therefore, $efg=n.$
\end{proposition}
In the paper \cite{Min_Sav}, Minculete and Savin introduced the notion of \textit{exponential divisor} of an ideal of the ring of integers of an algebraic number field.

Keeping the above notations, an \textit{exponential divisor} of the ideal $I=p\mathcal{O}_{K}$ has the following form
$d^{\left(e\right)}_{I}=P^{\beta_{1}}_{1}\cdot P^{\beta_{2}}_{2}\cdot...\cdot P^{\beta_{g}}_{g},$
where $\beta_{1},$ $\beta_{2},$...,$\beta_{g}\in\mathbb{N}^{*},$ with $\beta_{i}$ $|$ $e_{i},$ for $\left(\forall\right)$ $i=\overline{1,g}.$

 Let $\mathbb{J}$ be the set of ideals of the ring $\mathcal{O}_{K}.$ Minculete and Savin (in  \cite{Min_Sav}) extended the functions $\tau,$  $\tau^{\left(e\right)}$ to ideals of the ring $\mathcal{O}_{K}$ so:\\
$\tau : \mathbb{J}\rightarrow \mathbb{C},\tau\left(I\right)=$ the number of divisors of the ideal $I,$ respectively\\
$\tau^{\left(e\right)} : \mathbb{J}\rightarrow\mathbb{C}$, $\tau^{\left(e\right)}\left(I\right)=$ the number of exponential divisors of the ideal $I$. \\
Also, they gave in \cite{Min_Sav}, the formulas for $\tau\left(I\right)$  and for  $\tau^{\left(e\right)}\left(I\right),$ namely:

\begin{equation}\label{sec3_1}
\tau\left(I\right)=\left(e_{1}+1\right)\cdot\left(e_{2}+1\right)\cdot...\cdot\left(e_{g}+1\right)                       
\end{equation}

and 
\begin{equation}\label{sec3_2}
\tau^{\left(e\right)}\left(I\right)=\tau\left(e_{1}\right)\cdot \tau\left(e_{2}\right)\cdot...\cdot \tau\left(e_{g}\right)    
\end{equation}
(see Proposition 2.1 from \cite{Min_Sav}).
\smallskip

Let $K$ be an algebraic number field and let $\mathcal{O}_{K}$ be the ring of integers of $K$. Let $\mathbb{J}$ be the set of ideals of the ring $\mathcal{O}_{K}$ and let Spec$\left(\mathcal{O}_{K}\right)$ be the set of the prime ideals od the ring  $\mathcal{O}_{K}.$  Let $I$ be an ideal of the ring $\mathcal{O}_{K}.$  It exists and is unique $g\in\mathbb{N}^{*},$ the distinct ideals $P_{1},$ $P_{2},$..., $P_{g}\in Spec\left(\mathcal{O}_{K}\right)$ and the distincy numbers $e_{1},$ $e_{2},$...,$e_{g}\in\mathbb{N}^{*}$ such that $I=P^{e_{1}}_{1}\cdot P^{e_{2}}_{2}\cdot...\cdot P^{e_{g}}_{g}.$

We generalize  the notion of entropy of an ideal of like this:

\begin{definition}
\label{twodotone}
Let $I\neq \left(0\right)$ be an ideal of the ring  $\mathcal{O}_{K}$, decomposed as above. We define the entropy of the ideal $I$ as follows:
\begin{equation}\label{sec3_3}
H\left(I\right)= - \sum_{i=1}^{g}\frac{e_{i}}{\Omega(I)} \log \: \frac{e_{i}}{\Omega(I)},     
\end{equation}
where $ \Omega\left(I\right) =  e_{1}+ e_{2}+...+ e_{g}.$
\end{definition}
Immediately, we obtain the following equivalent form, for the entropy of the ideal $I:$
\begin{equation}\label{sec3_4}
H\left(I\right)=\log \: \Omega\left(I\right) - \frac{1}{\Omega\left(I\right)}\cdot \sum_{i=1}^{g}e_{i}\cdot \log \: e_{i}.  
\end{equation}
We now give some examples of calculating the entropy of an ideal.
\begin{example}
\label{twodottwo}
Let $\xi$ be a primitive root of order $5$ of the unity and let $K=\mathbb{Q}\left(\xi\right)$ be the $5$th cyclotomic field. The ring of algebraic integers of the field $K$ is $\mathcal{O}_{K}=\mathbb{Z}\left[\xi\right].$
We consider the ideal $\left(1-\xi\right)\cdot \mathbb{Z}\left[\xi\right].$ It is known that $\left(1-\xi\right)\cdot \mathbb{Z}\left[\xi\right]\in Spec\left(\mathcal{O}_{K}\right)$ (see \cite{ireland},  \cite{Savin}). 
Let the ideal $5\cdot \mathbb{Z}\left[\xi\right]=\left(1-\xi\right)^{4}\cdot \mathbb{Z}\left[\xi\right].$ The entropy of the ideal $5\cdot \mathbb{Z}\left[\xi\right]$ is 
$$H\left(5\cdot \mathbb{Z}\left[\xi\right]\right)=\log \: 4- \frac{1}{4}\cdot 4 \cdot \log \: 4=0.$$
\end{example}
\begin{example}
\label{twodotthree}
Let the pure cubic field $K=\mathbb{Q}\left(\sqrt[3]{2}\right).$ Since $2^{2}\not\equiv 1$ (mod $9$), it results that the ring of algebraic integers of the field $K$ is  $\mathcal{O}_{K}=\mathbb{Z}\left[\sqrt[3]{2}\right]$
(see \cite{ Murty}).\\
Since $29\equiv 2$ (mod $3$), then $29\mathbb{Z}\left[\sqrt[3]{2}\right]=P_{1}\cdot P_{2},$ where $P_{1}, P_{2}\in Spec \left(\mathbb{Z}\left[\sqrt[3]{2}\right]\right).$ So, the ideal $29\mathbb{Z}\left[\sqrt[3]{2}\right]$ splits in the ring $\mathbb{Z}\left[\sqrt[3]{2}\right].$
The entropy of the ideal $29\mathbb{Z}\left[\sqrt[3]{2}\right]$ is 
$$H\left(29\mathbb{Z}\left[\sqrt[3]{2}\right]\right)=\log \: 2 - \frac{1}{2}\cdot 2 \cdot \log \: 1=\log \: 2.$$

\end{example}

\begin{example}
\label{twodotfour}
In the same field (as in the previous example) $K=\mathbb{Q}\left(\sqrt[3]{2}\right)$ with the ring of integer $\mathcal{O}_{K}=\mathbb{Z}\left[\sqrt[3]{2}\right],$ we consider the ideal $31\mathbb{Z}\left[\sqrt[3]{2}\right].$\\
Since $31\equiv 1$ (mod $3$), then $31\mathbb{Z}\left[\sqrt[3]{2}\right]=P_{1}\cdot P_{2}\cdot P_{3},$ where $P_{1}, P_{2}, P_{3}\in Spec\left(\mathbb{Z}\left[\sqrt[3]{2}\right]\right).$ So, the ideal $31\mathbb{Z}\left[\sqrt[3]{2}\right]$ splits completely in the ring $\mathbb{Z}\left[\sqrt[3]{2}\right]$
(see \cite{ Murty}). The entropy of the ideal $31\mathbb{Z}\left[\sqrt[3]{2}\right]$ is 
$$H\left(31\mathbb{Z}\left[\sqrt[3]{2}\right]\right)=\log \: 3 - \frac{1}{3}\cdot 3 \cdot \log \: 1=\log \: 3.$$

\end{example}
\begin{remark}
\label{twodotfour}
Let $K$ be an algebraic number field and  let $\mathcal{O}_{K}$ be its  ring of integers.  Let $p$ be a prime positive integer. If $p$  is inert or totally ramified in the ring  $\mathcal{O}_{K},$  then $H\left(p\mathcal{O}_{K}\right)=0.$

\end{remark}

\begin{proof}
To calculate the entropy of ideal $p\mathcal{O}_{K}$, we consider two cases.\\
\textbf{Case 1}: if $p$ is inert in the ring  $\mathcal{O}_{K},$  it results that $p\mathcal{O}_{K}$ is a prime ideal. Then $ \Omega\left(p\mathcal{O}_{K}\right) =  1$ and $H\left(p\mathcal{O}_{K}\right)=0.$\\
\textbf{Case 2}: if $p$ is totally ramified in the ring  $\mathcal{O}_{K},$  it results that $p\mathcal{O}_{K}= P^{n},$ where $P\in Spec \left(\mathcal{O}_{K}\right)$ and $n=[K :\mathbb{Q}].$
It results immediately $ \Omega\left(p\mathcal{O}_{K}\right) =  n$ and $H\left(p\mathcal{O}_{K}\right)=log \:n - log \:n=0.$
\end{proof}
\begin{proposition}
\label{twodotfive}
 Let $n$ be a positive integer, $n\geq2$ and let $p$ be a positive prime integer. Let $K$ be an algebraic number field of degree $[K :
\mathbb{Q}] = n$ and let $\mathcal{O}_{K}$ be its the ring of integers. Then:\\
 \begin{equation}\label{sec3_5}
0\leq H\left(p\mathcal{O}_{K}\right)\leq \log \: \omega\left(p\mathcal{O}_{K}\right)\leq  \log \:n, \; 
\end{equation}
where $ \omega\left(p\mathcal{O}_{K}\right)$  is the number of distinct prime factors of the ideal $p\mathcal{O}_{K}.$
\end{proposition}
\begin{proof}
The proof of the inequality $0\leq H\left(p\mathcal{O}_{K}\right)\leq \log \: \omega\left(p\mathcal{O}_{K}\right)$ is similar to the proof of  Proposition \ref{onedotone} (that is, Theorem 2. from the article \cite{Minculete1}).\\
Since $\mathcal{O}_{K}$ is a Dedekind ring, the ideal $p\mathcal{O}_{K}$ is written in a unique way:
$$p\mathcal{O}_{K}=P^{e_{1}}_{1}\cdot P^{e_{2}}_{2}\cdot...\cdot P^{e_{g}}_{g},$$
where  $g\in\mathbb{N}^*,$   $e_{1}, e_{2},..., e_{g}\in\mathbb{N}^*$  and $P_{1},$ $P_{2},$..., $P_{g}\in Spec\left(\mathcal{O}_{K}\right).$ Applying Proposition  \ref{onedotthree} i) we obtain that $\omega\left(p\mathcal{O}_{K}\right)=g\leq n.$
The equality $\omega\left(p\mathcal{O}_{K}\right)=n$ is achieved when the ideal $p$ splits totally in the ring $\mathcal{O}_{K}.$ It follows that 
 $$0\leq H\left(p\mathcal{O}_{K}\right)\leq log \: \omega\left(p\mathcal{O}_{K}\right)\leq  log \:n.$$
 
\end{proof}
\begin{proposition}
\label{twodotsix}
Let $K$ be an algebraic number field and let $\mathcal{O}_{K}$ be its the ring of integers. Let $p$ be a prime positive integer. If the extension of fields $\mathbb{Q}\subset K$ is a Galois extension, then
$$H\left(p\mathcal{O}_{K}\right)=\log \: \omega\left(p\mathcal{O}_{K}\right).$$
\end{proposition}
\begin{proof}
Taking into account the fact that  $\mathcal{O}_{K}$ is a Dedekind rimg and applying Proposition \ref{onedotthree} ii), it follows that the ideal $p\mathcal{O}_{K}$ is uniquely written as follows:

$$p\mathcal{O}_{K}=P^{e_{1}}_{1}\cdot P^{e_{1}}_{2}\cdot...\cdot P^{e_{1}}_{g},$$
where  $g\in\mathbb{N}^*,$   $e_{1}\in\mathbb{N}^*$  and $P_{1},$ $P_{2},$..., $P_{g}\in Spec\left(\mathcal{O}_{K}\right).$ According to formula (1.2), the entropy of the ideal $p\mathcal{O}_{K}$ is
$$H\left(p\mathcal{O}_{K}\right)=\log\left(ge_{1}\right) - \frac{1}{ge_{1}}\cdot ge_{1} \cdot \log \: e_{1}=  \log \: g= \log \: \omega\left(p\mathcal{O}_{K}\right). $$

\end{proof}

\smallskip

\section{The entropy and exponential divisiors}

\noindent 
\smallskip\\
\begin{proposition}
\label{fourdotone}
Let $\alpha,\beta$ be two natural numbers with $\alpha\geq\beta\geq 1$, $p$ a prime number and decomposition in prime factors of $n$ given by $n=\prod_{i=1}^rp_i^{a_i}$ with $a_i\in\mathbb{N}^*$ for all $i\in\{1,...,r\}$ such that $\gcd(n,p)=1$. Then\\
i) If $\beta\geq\Omega(n)e^{-H(n)}$, then we have $H(np^\alpha)\leq H(np^\beta)$;\\
ii) If $\beta\leq\alpha\leq\Omega(n)e^{-H(n)}$, then we have $H(np^\alpha)\geq H(np^\beta)$;\\
iii) If $\beta\leq\Omega(n)e^{-H(n)}\leq\alpha$, then we have $H(np^\alpha)\leq H(np^\beta)$.
\end{proposition}
\begin{proof}
We study the difference of the entropy of the numbers $np^\alpha$ and $np^\alpha$, with $\gcd(n,p)=1$ and $\alpha\geq\beta\geq 1$. Using equality \eqref{ms_5}, we deduce
\begin{align*}
H(np^\alpha)-H(np^\beta)&=\frac{(\alpha-\beta)(\Omega(n)\log\Omega(n)-\Omega(n)H(n))}{(\Omega(n)+\alpha)(\Omega(n)+\beta)}+\log\frac{\Omega(n)+\alpha}{\Omega(n)+\beta}\\
&+\frac{\beta\log\beta}{\Omega(n)+\beta}-\frac{\alpha\log\alpha}{\Omega(n)+\alpha}\\
&=\frac{(\alpha-\beta)\sum_{i=1}^{r}\alpha_{i}\log \alpha_{i}}{(\Omega(n)+\alpha)(\Omega(n)+\beta)}+\log\frac{\Omega(n)+\alpha}{\Omega(n)+\beta}\\
&+\frac{\beta\log\beta}{\Omega(n)+\beta}-\frac{\alpha\log\alpha}{\Omega(n)+\alpha}.
\end{align*}
We note $\sum_{i=1}^{r}\alpha_{i}\log \alpha_{i}=A\geq 0$ and $\Omega(n)=t\geq 1$ and we take the function $f:[\beta,\infty)\to\mathbb{R}$ defined by $f(\alpha)=\frac{(\alpha-\beta)A}{(t+\alpha)(t+\beta)}+\log\frac{t+\alpha}{t+\beta}+\frac{\beta\log\beta}{t+\beta}-\frac{\alpha\log\alpha}{t+\alpha}.$ But, $f'(\beta)=\frac{A-t\log\alpha}{(t+\alpha)^2}=\frac{t}{(t+\alpha)^2}(\frac{A}{t}-\log\alpha)=\frac{\Omega(n)}{(\Omega(n)+\alpha)^2}(\log\Omega(n)-H(n)-\log\alpha),$ which is equivalent to $f'(\alpha)=\frac{\Omega(n)}{(\Omega(n)+\alpha)^2}\log\frac{\Omega(n)e^{-H(n)}}{\alpha}.$
In the first case $\beta\geq\Omega(n)e^{-H(n)}$, then we deduce $f'(\alpha)\leq 0$, so, the function $f$ is decreasing. Therefore for $\alpha\geq\beta$ we have $f(\alpha)\leq f(\beta)=0$, which means that $H(np^\alpha)\leq H(np^\beta)$. The second case $\beta\leq\alpha\leq\Omega(n)e^{-H(n)}$, then $f'(\alpha)\geq 0$, so, the function $f$ is increasing. Since $\alpha\geq\beta$ we deduce $f(\alpha)\leq f(\beta)=0$, which means that $H(np^\alpha)\geq H(np^\beta)$. The third case is similar to the first case. Consequently the statement is true.
\end{proof}
\begin{remark}
\label{remark}
 Let $n$ be a positive integer,  $n\geq2,$ $n=p^{\alpha_{1}}_{1}p^{\alpha_{2}}_{2}...p^{\alpha_{r}}_{r}$, where $p_{1}, p_{2},...,p_{r}$ are distinct prime positive integers and $\alpha_{1}, \alpha_{2},...,\alpha_{r}\in\mathbb{N}^{*}.$ 
 If $n$ is a square free integer, then $\alpha_{1}= \alpha_{2}=...=\alpha_{r}=1.$ So, $\beta\geq 1=\Omega(n)e^{-H(n)}$, then we have $H(np^\alpha)\leq H(np^\beta)$. Therefore, for any $d_{e}$ be an exponential divisor of $m=np^{\alpha}$, the following inequality hold:
$$H\left(d_{e}\right)\geq H\left(m\right).$$
If $n$ is not a square free integer, then for an exponential divisor $d_{e}$ of $m=np^{\alpha}$ given by $d_{e}=np^{\beta}$ , the following inequality hold:
$$H\left(d_{e}\right)\leq H\left(m\right),$$
when $\beta\leq\alpha\leq\Omega(n)e^{-H(n)}$ or $$H\left(d_{e}\right)\geq H\left(m\right),$$
when $\beta\leq\Omega(n)e^{-H(n)}\leq\alpha.$
\end{remark}

\smallskip

Considering now a positive integer $n$ of the form  $n=p^{\alpha_{1}}_{1}\cdot p^{\alpha_{2}}_{2}\cdot...\cdot p^{\alpha_{r}}_{r}$, where $p_{1}, p_{2},...,p_{r}$ are distinct prime positive integers and $\alpha_{1}, \alpha_{2},...,\alpha_{r}\in$ $\left\{1, 2\right\}$ and an arbitrary exponential divisor $d_{e}$ of $n$, we are interested in the relationship between the entropy $ H\left(n\right)$ and the entropy $H\left(n\right).$\\
For example, when $r=2$ and $n=12=2^{2}\cdot 3,$ $ H\left(12\right)=\log \: 3 -  \frac{2}{3}\cdot \log \: 2.$ $d_{e}=6$ is an exponential divisor of $12$ and  $ H\left(6\right)=  \log \: 2,$ so $ H\left(12\right)\leq H\left(6\right)$ (we are in the case iii) from Proposition \ref{fourdotone}).\\
Another example: when $r=3$ and $n=180=2^{2}\cdot 3^{2}\cdot 5,$ $ H\left(180\right)=\log \: 5 -  \frac{4}{5}\cdot \log \: 2.$ $d_{e}=60=2^{2}\cdot 3\cdot 5$ is an exponential divisor of $180$ and  $ H\left(60\right)=  \log \: 4 -  \frac{1}{2}\cdot \log \: 2,$ so $ H\left(60\right)\leq H\left(180\right)$ (we are in the case ii). from Proposition \ref{fourdotone}).\\

\begin{corollary}
\label{foudottwo}
 Let $n$ be a positive integer,  $n\geq2,$ $n=p^{\alpha_{1}}_{1}\cdot p^{\alpha_{2}}_{2}\cdot...\cdot p^{\alpha_{r}}_{r}$, where $r\geq3,$ $p_{1}, p_{2},...,p_{r}$ are distinct prime positive integers and $\alpha_{1}, \alpha_{2},...,\alpha_{r}\in$ $\left\{1, 2\right\}.$ Then, for any $d_{e}$ be an exponential divisor of $n,$ the following inequality hold:
$$H\left(d_{e}\right)\leq H\left(n\right).$$
\end{corollary}
\begin{proof}
Let $d_{e}$ be an exponential divisor of $n,$ then $d_{e}=p^{\beta_{1}}_{1}\cdot p^{\beta_{2}}_{2}\cdot...\cdot p^{\beta_{r}}_{r}$, where  $\beta_{1}, \beta_{2},$
$...,\beta_{r}\in\mathbb{N}^{*},$  with $ \beta_{i}\in\left\{1, 2\right\},$  $ \beta_{i} |  \alpha_{i}$ $\left(\forall\right)$ $i=\overline{1,r}.$ 
Without diminishing the generality, we assume that there exists $s\in\mathbb{N}^{*},$ $s\leq r$ such that $\alpha_{1}=...=\alpha_{s}=2$ and $\alpha_{s+1}=...=\alpha_{r}=1$ and 
 there exists $l\in\mathbb{N}^{*},$ $l\leq s$ such that $\beta_{1}=...=\beta_{l}=2$ and $\beta_{l+1}=...=\beta_{r}=1.$ Then $n=p^{2}_{1}\cdot p^{2}_{2}\cdot...\cdot p^{2}_{l}\cdot....\cdot p^{2}_{s}\cdot p_{s+1}\cdot...\cdot p_{r}$ and  $d_{e}= p^{2}_{1}\cdot p^{2}_{2}\cdot...\cdot p^{2}_{l}\cdot  p_{l+1}\cdot  ...\cdot p_{s} \cdot  p_{s+1}\cdot...\cdot p_{r}.$
Repeatedly applying Proposition \ref{fourdotone} or Remark \ref{remark}, we have:

$$ H\left(n\right)=H\left(p^{2}_{1}\cdot p^{2}_{2}\cdot...\cdot p^{2}_{l}\cdot....\cdot p^{2}_{s-1}\cdot p^{2}_{s}\cdot p_{s+1}\cdot...\cdot p_{r}\right)\leq $$
$$\leq H\left(p^{2}_{1}\cdot p^{2}_{2}\cdot...\cdot p^{2}_{l}\cdot....\cdot p^{2}_{s-1}\cdot p_{s}\cdot p_{s+1}\cdot...\cdot p_{r}\right)\leq ..........\leq  $$
$$\leq H\left(p^{2}_{1}\cdot p^{2}_{2}\cdot...\cdot p^{2}_{l}\cdot  p_{l+1}\cdot...\cdot p_{s}\cdot p_{s+1}\cdot...\cdot p_{r}\right)=H\left(d_{e}\right).     $$

\end{proof}

\begin{corollary}
\label{foudotthree}
Let $K$ be an algebraic number field and let $\mathcal{O}_{K}$ be the ring of integers of $K.$  Let $I$ be an ideal of the ring $\mathcal{O}_{K},$  $I=P^{e_{1}}_{1}\cdot P^{e_{2}}_{2}\cdot...\cdot P^{e_{g}}_{g},$ where $g\geq3,$ $P_{1},$ $P_{2},$..., $P_{g}$ are distinct prime ideals of the ring $\mathcal{O}_{K}$ 
and $e_{1}, e_{2},...,e_{r}\in\left\{1, 2\right\}.$  Then, for any  $d^{\left(e\right)}_{I}$ be an exponential divisor of the ideal $I,$ the following inequality hold:
$$    H\left(d^{\left(e\right)}_{I}\right) \leq H\left(I\right).$$

\end{corollary}
\begin{proof}
The proof is similar to the proof of Corollary  \ref{foudottwo}.

\end{proof}

\noindent 
\smallskip\\


\begin{thebibliography}{99}





\smallskip
 
\smallskip




\bibitem{ireland}
K. Ireland,  M. Rosen,  \textit{A Classical Introduction to Modern Number Theory}, Springer Verlag, 1992.

\bibitem{Jeong}
S. Jeong, K. H. Kim and G. Kim,  \textit{Algebraic entropies of natural numbers with one or two factors}, J. Korean Soc. Math. Educ. Ser. B: Pure Appl. Math., Vol. 23. No. 3. 2016, p.205-221.

\bibitem{Minculete1}
N. Minculete,   C. Pozna, \textit{The Entropy of a Natural Number}, Acta Technica Jaurinensis, Vol. 4. No. 4. 2011, p.425-431.

\bibitem{Minculete}
N. Minculete,  \textit{On  certain  inequalities  about  arithmetic  functions  which  use  the 
exponential  divisors}, Int. J. Number  Theory, 8, Issue 6, 2012, 1527-1535. 

\bibitem{Min_Sav} 
N. Minculete, D. Savin,  \textit{Some generalizations of the functions} $\tau$ \textit{and} $\tau^{(e)}$ \textit{in algebraic number fields},  Expositiones Mathematicae, vol. 39 (2021), p. 344 -353.

\bibitem{Murty}
M. R. Murty, J. Esmonde, \textit{Problems in algebraic number theory}, Second Edition, Springer, 2005.

\bibitem{Pet}
Y.-F.S. P\'etermann and J. Wu, \textit{On the sum of exponential divisors of an integer}, Acta Math. Hungar. 77 (1997), 159-175. 

\bibitem{Savin}
D. Savin, M. \c{S}tefanescu, \textit{Lessons of Arithmetics and Number Theory}, Matrix Rom Publishing House: Bucharest, Romania, 2008 (In Romanian).

\bibitem{Sandor_2}
J. S\'andor and B. Crstici, \textit{Handbook of Number Theory II}, Kluwer Academic Publishers, Dordrecht/Boston/London, 2004.

\bibitem{Sandor}
J. S\'andor, \textit{A Note on Exponential Divisors and Related Arithmetic Functions}, Scientia Magna, Vol 1 (2006), no. 1.

\bibitem{Straus} E. G. Straus and M. V. Subbarao, \textit{On exponential divisors}, Duke Math. J. 41 (1974) 465-471.

\bibitem{Subbarao}
M. V. Subbarao, \textit{On some arithmetic convolutions in The Theory of Arithmetic Functions}, Lecture Notes in Mathematics, New York, Springer-Verlag, 1972.

\bibitem{Toth}
L. T\'oth, \textit{On Certain Arithmetic Functions Involving Exponential Divisors}, Annales Univ. Sci. Budapest., Sect. Comp. 24 (2004), 285-294.

\bibitem{Wu} 
J. Wu, \textit{Probl\`eme de diviseurs exponentiels et entiers exponentiellement sans facteur carr\'e}, Journal de Th\'eorie des Nombres de Bordeaux, tome 7, no 1 (1995), p. 133-141.

\bibitem{link}
https://math.stackexchange.com/questions/2369779/entropy-of-a-natural-number
\end{thebibliography}
\end{document}